\documentclass[10pt]{article} 
\usepackage{latexsym}
\usepackage{amsfonts}
\usepackage{amssymb}
\usepackage{graphicx}

\begin{document}

\title{MATRIX REPRESENTATIONS OF OCTONIONS AND THEIR  APPLICATIONS}
\author{Yongge   Tian \\
Department of Mathematics and Statistics \\
Queen's University \\
Kingston,  Ontario,  Canada K7L 3N6\\
{\tt e-mail:ytian@mast.queensu.ca}}
\date{}
\maketitle

\begin{abstract} As is well-known, the real quaternion division
algebra $\Bbb H$ is algebraically isomorphic to a 4-by-4 real matrix
algebra. But the real division octonion algebra $\Bbb O$ can not be
algebraically isomorphic to any matrix algebras over the real number field
$\Bbb R$, because $\Bbb O$ is a non-associative algebra over
 $\Bbb R$. However since $\Bbb O$ is an extension of $\Bbb H$
 by the Cayley-Dickson process and is also finite-dimensional, some
 pseudo real matrix representations of octonions can still be
 introduced through real matrix representations of quaternions. In this
 paper we give a complete investigation to real matrix representations of
 octonions, and consider their various applications to octonions as well as 
 matrices of octonions. \\

\end{abstract}
\medskip

\noindent {\Large {\bf 1. \ Introduction }} \\

\noindent Let $ \Bbb O$ be the octonion algebra over the real number
 field $\Bbb R$. Then it is well known by the Cayley-Dickson process that any
 $ a \in \Bbb O$ can be written as
$$ 
a = a' + a''e, \eqno (1.1)
$$
where $a', \,  a'' \in  \Bbb H = \{ \, a = a_0 + a_1i + a_2j + a_3k \ |
 \ i^2 = j^2 = k^2 = -1, \ ijk = -1, \ a_0$---$a_3 \in \Bbb R \, \}$,
 the real quaternion division algebra. The addition and multiplication
 for any $ a = a' + a''e, \  b = b' + b''e \in \Bbb O$ are defined by
$$ 
a  + b = ( \, a' + a''e \,)  + ( \, b' + b''e ) =  (  \, a' + b' \,) +
( \, a'' +  b'' \, )e,
 \eqno(1.2)
$$  
and
$$
ab = ( \,a' + a''e )(b' + b''e) = ( a'b'  -   \overline{b''}a'' ) +
(  b''a' +  a''\overline{b'} \, )e,   \eqno(1.3)
$$ 
where $ \overline{b'}, \, \overline{b''}$ denote the conjugates  of the
quaternions $ b'$ and $ b''$.  In that case, $\Bbb O$  is an
eight-dimensional non-associative but alternative division algebra over
its center field $\Bbb R $, and  the canonical basis of $\Bbb O$ is 
$$
 1, \ \ \  e_1 =i, \ \ \  e_2 = j, \ \ \  e_3 = k, \ \ \  e_4 =e, \ \ \
 e_5 = ie, \ \  e_6 = je, \ \ \  e_7 = ke.    \eqno(1.4)
$$
The multiplication rules for the basis of $\Bbb O$ are listed in the
following matrix 
$$
E_8^TE_8 = \left[ \begin{array}{crrrrrrr}  1  & e_1 &  e_2 &  e_3 &  e_4 &
 e_5 &  e_6 &  e_7
 \\   e_1 & -1 &  e_3 &  -e_2 &  e_5 &  -e_4 &  -e_7 & e_6 
 \\ e_2 &  -e_3 &  -1 &  e_1 &  e_6 &  e_7 &  -e_4 & -e_5  \\
  e_3 &  e_2 &  -e_1 & -1  &   e_7 &  -e_6 &  e_5 & -e_4  \\   e_4
   & -e_5 &  -e_6 &  -e_7 &  -1 &  e_1 &  e_2 &  e_3  \\   e_5  & e_4 &
   -e_7 &  e_6 & -e_1 &  -1 &  -e_3 &  e_2  \\
   e_6  & e_7 &  e_4 &  -e_5 &  -e_2 &  e_3 &  -1 &  -e_1  \\ e_7  & -e_6 &
   e_5 &  e_4 &  -e_3 &  -e_2 &  e_1 &  -1  \end{array} \right],
  \eqno(1.5)
$$ 
where $ E_8 = [ \, 1,  \ e_1, \ \cdots, \  e_7 \,]$.  Under Eq.(1.4)
all elements of $\Bbb O$ take the form
$$  
a = a_0 + a_1e_1 + \cdots + a_7e_7,  \eqno(1.6)
$$ 
where $ a_0$---$a_7 \in \Bbb R$, which can also simply be written as
$a = {\rm Re\,}a +  {\rm Im\,}a,$
where $ {\rm Re\,}a = a_0.$ The conjugate of $ a $ is defined to be 
 $$ 
\overline{a} = \overline{a'} - a''e = {\rm Re\,}a -  {\rm Im\,}a.
\eqno(1.7)
$$
This operation satisfies
$$ 
\overline{\overline{a}} = a,  \ \ \  \overline{a+b} =  \overline{a} +
\overline{b},  \ \ \  \overline{ab} =   \overline{b}\overline{a}
\eqno(1.8)
$$
for all $ a, \, b  \in \Bbb O.$   The norm of $ a $ is defined to be
 $ |a | := \sqrt{a\overline{a}} = \sqrt{\overline{a}a} =
 \sqrt{ a_0^2 + a_1^2 + \cdots + a_7^2 }.$
Although $ \Bbb O $ is nonassociative, it is still an alternative,
flexible, quadratic, composition and division algebra over $\Bbb R,$
that is, for all $ a, \, b \in \Bbb O$, the following equalities hold:
$$
\displaylines{
\hspace*{2cm}
a(ab) = a^2b,  \qquad (ba)a = ba^2,   \qquad (ab)a = a(ba) := aba,
\hfill (1.9)
\cr
\hspace*{2cm}
a^{-1} = \frac{ \overline{a}}{ |a|^2},    \hfill (1.10)
\cr
\hspace*{2cm}
a^2 - 2({\rm Re\,}a)a + |a|^2 = 0,  \qquad  ({\rm Im\,}a)^2 = -|
{\rm Im\,}a|^2,  \hfill (1.11)
\cr
\hspace*{2cm}
|ab| = |a||b|.
\hfill (1.12)
\cr}
$$

As is well known, any finite-dimensional associative algebra over an
arbitrary field $\Bbb F$ is algebraically isomorphic to a subalgebra of a
total matrix algebra over the field. In other words, any element in a
 finite-dimensional associative algebra over $\Bbb F $ has a faithful matrix
 representation over the field. For the real quaternion algebra
 $\Bbb H$, it is well known that
 through the bijective map
$$
\phi : a= a_0 + a_1i + a_2j + a_3k  \in \Bbb H \longrightarrow
 \phi(a) = \left[ \begin{array}{rrrr} a_0 & -a_1 &  -a_2 & - a_3  \\ a_1 & a_0 & - a_3 & a_2 \\
a_2 & a_3 &  a_0 &  -a_1  \\ a_3 & -a_2 &  a_1 & a_0  \end{array} \right],
\eqno (1.13)
$$
$\Bbb H$ is algebraically isomorphic to the matrix algebra
$$
{\cal M} = \left\{ \, \left. \left[ \begin{array}{rrrr} a_0 & -a_1 &  -a_2 & - a_3  
\\ a_1 & a_0 & - a_3 & a_2 \\
a_2 & a_3 &  a_0 &  -a_1  \\ a_3 & -a_2 &  a_1 & a_0  \end{array} \right]
 \, \right|\, a_0, \ a_1, \ a_2, \ a_3 \in \Bbb R  \right\},
\eqno (1.14)
$$
and $\phi(a)$ is a faithful real matrix representation of $a$.
Our consideration for matrix representations of octonions are
based on Eqs.(1.1)---(1.3) and the result in Eq.(1.13).

\medskip


We next present some basic results related to matrix representations of
quaternions, which will be serve as a tool for our examination in the
sequel.

\medskip

\noindent {\bf Lemma 1.1}\cite{Ti}. \, {\em Let $ a= a_0 + a_1i + a_2j + a_3k \in
\Bbb H$ be given$,$ where $ a_0 $---$ a_3
\in \Bbb R $.  Then the diagonal matrix  $ {\rm diag }( \, a, \ a, \ a, \
a \, )$  satisfies the following unitary similarity factorization equality
$$
 Q \left[ \begin{array}{cccc} a &  &  &  \\  & a &  & \\
 &  &  a &    \\  & & & a  \end{array} \right] Q^* =
 \left[ \begin{array}{rrrr} a_0 & -a_1 &  -a_2 & - a_3  \\ a_1 & a_0 &
 - a_3 & a_2 \\
a_2 & a_3 &  a_0 &  -a_1  \\ a_3 & -a_2 &  a_1 & a_0  \end{array}
\right] \in \Bbb R^{4 \times 4},
 \eqno (1.15)
$$ 
where the matrix $ Q $ has the following independent expression   
$$
 Q = Q^* = \frac{1}{2} \left[ \begin{array}{rrrr} 1 & i & j & k \\ -i & 1 & k & -j \\ -j & -k & 1 &  i  \\  -k & j & -i & 1  \end{array} \right], \eqno (1.16)
$$  
which is a unitary matrix over $\Bbb H$.} 

\medskip

\noindent {\bf Lemma 1.2}\cite{Ti}. \, {\em  Let $a, \,  b \in \Bbb H$, and
$\lambda \in \Bbb R$. Then

{\rm (a)} \ $ a=b  \Longleftrightarrow  \phi(a) = \phi(b).$

{\rm (b)}  \ $ \phi(a+b)= \phi(a) + \phi(b), \ \ \ \phi(ab)= \phi(a) \phi(b),
\ \ \ \phi( \lambda a)=  \lambda \phi(a),  \ \ \ \phi(1) = I_4.$ 

{\rm (c)}  \ $ a =  \frac{1}{4} E_4  \phi(a)E_4^*,$ where
$ E_4 := [ \, 1, \ i, \  j, \  k \, ]$ and $E_4^* := [ \, 1, \ -i, \  -j, \
 -k \, ]^T. $

{\rm (d)}  \ $  \phi(\overline{a})= \phi^T(a). $ 

{\rm (e)} \  $ \phi(a^{-1})= \phi^{-1}(a),$  if $  a \neq 0.$ 

{\rm (f)}  \  $ {\rm det\,} [\phi(a)]= |a|^4 . $   }

\medskip

We can also introduce from Eq.(1.13) another real matrix representation
of $ a $ as follows
$$ 
\tau(a) :=  K\phi^T(a)K =  \left[ \begin{array}{rrrr} a_0 & -a_1 &  -a_2 &
 -a_3  \\ a_1 & a_0 & a_3 & -a_2 \\  a_2 & -a_3 &  a_0 &  a_1  \\ a_3 & a_2
 &  -a_1 & a_0  \end{array} \right],  \eqno (1.17)
$$
where $ K = {\rm diag}( \, 1, \ -1, \ -1, \ -1 \, )$. Some  basic operation
properties on $ \tau(a)$  are
$$
  \tau(a+b) =\tau(a) + \tau(b), \qquad \tau(ab) =\tau(b)\tau(a), \qquad
  \tau(\overline{a}) = \tau^T(a), \eqno (1.18)
 $$
$$
 {\rm \det\,}[\phi(a)]= |a|^4 , \qquad  \phi(a^{-1})= \phi^{-1}(a) \ \ \
{\rm  if } \ \  a \neq 0.   \eqno (1.19)
$$  

Combining the two real matrix representations of quaternions with their real
 vector representations, we have the following important result. 

\medskip

\noindent {\bf Lemma 1.3.} \, {\em  Let $  x = x_0 + x_1i + x_2j + x_3k
\in  \Bbb H,$  and denote $ \overrightarrow{x}
 = [\, x_0, \ x_1, \  x_2, \  x_3 \, ]^T,$ called the vector representation of $ x $. Then  for  all 
$a, \, b, \, x \in \Bbb H$, we have 
$$ 
\overrightarrow{ax} = \phi(a) \overrightarrow{x}, \qquad \overrightarrow{xb}
= \tau(b) \overrightarrow{x},  \qquad \overrightarrow{axb} = \phi(a) \tau(b)
\overrightarrow{x} = \tau(b) \phi(a) \overrightarrow{x},
 \eqno (1.20) 
$$
and the equality 
$$ 
\phi(a) \tau(b) = \tau(b) \phi(a)   \eqno (1.21) 
$$ 
always holds.} 

\medskip

\noindent {\bf Proof.}\, Observe that 
$$ 
 \overrightarrow{x}= \phi(x) \alpha^T_4, \qquad  \overrightarrow{x}= \tau(x) \alpha^T_4,  \ \ \ \  \alpha_4 = [ \, 1, \ 0, \ 0, \ 0 \, ].
$$
We find by Lemma 1.1 and Eq.(1.2) that
$$
 \overrightarrow{ax} = \phi(ax) \alpha^T_4 = \phi(a) \phi(x) \alpha^T_4
  = \phi(a) \overrightarrow{x}, \qquad  \overrightarrow{xb} = \tau(xb)
  \alpha^T_4 = \tau(b) \tau(x) \alpha^T_4
  = \tau(b) \overrightarrow{x},  
$$ 
and
$$
\overrightarrow{axb} =  \overrightarrow{a(xb)} = \phi(a) \overrightarrow
{(xb)} = \phi(a) \tau(b) \overrightarrow{x},  \ \ \ \overrightarrow{axb} =
\overrightarrow{(ax)b} = \tau(b) \overrightarrow{(ax)} = \tau(b)\phi(a)
\overrightarrow{x}.
$$ 
These four equalities are exactly the results in Eqs.(1.20) and (1.21).
\qquad $ \Box $ 

\medskip









\noindent {\bf Lemma 1.4}\cite{Sc}\cite{Zhe}. \, {\em  Let $a, \, b, \, x  \in \Bbb O$ be given. Then

{\rm (a)} \  $  {\rm Re\,}(ab) = {\rm Re\,}(ba), \qquad  {\rm Re\,}((ax)b)
= {\rm Re\,}(a(xb)).$

{\rm (b)} \  $ (aba)x  = a(b(ax)), \qquad x(aba)  = ((xa)b)a.$ 

{\rm (c)} \  $ (ab)(xa) = a(bx)a, \qquad (bx)(ab) = b(xa)b.$ 

{\rm (d)} \ $ ( a, \  b, \ x ) = - ( a, \  x, \ b ) = ( x, \ a, \ b ),$   where $ ( a, \ b, \  x ) = 
(ab)x - a(bx). $ } \\

\noindent  {\Large {\bf 2. \ The real matrix representations of octonions }} \\

\noindent Based on the results on the real matrix representation of
quaternions, we now can introduce  real matrix representation of
octonions.
  
\medskip

\noindent {\bf Definition 2.1.} \,  Let $ a = a' +  a''e \in \Bbb O,$
where $ a' = a_0 + a_1i + a_2j + a_3k, \, a'' = a_4 + a_5i +  a_6j + a_7k \in
\Bbb H.$ Then the $ 8 \times 8 $ real matrix
$$ 
\omega(a) :=  \left[ \begin{array}{cc} \phi(a') &  - \tau(a'') K_4   \\  
 \phi(a'')K_4  & \tau(a')
  \end{array} \right],  \eqno (2.1)
$$     
is called the left matrix representation of $ a $ over $\Bbb R,$ where $ K_4 = {\rm diag}( 1, \ -1, \ -1, \
 -1 )$. Written  in an explicit form, 
$$ 
\omega(a) =  \left[ \begin{array}{crrrrrrr}  
a_0  & -a_1 &  -a_2 &  -a_3 &  -a_4 &  -a_5 &  -a_6 &  -a_7 \\   
a_1 & a_0 &  -a_3 &  a_2 &  -a_5 &  a_4 &  a_7 & -a_6  \\ 
a_2 &  a_3 & a_0 &  -a_1 &  -a_6 &  -a_7 &  a_4 & a_5  \\
  a_3 &  -a_2 &  a_1 & a_0  &   -a_7 &  a_6 &  -a_5 & a_4  \\   a_4  & a_5 &  a_6 &  a_7 &  a_0 &  -a_1 &  -a_2 &  -a_3  \\   a_5  & -a_4 &  a_7 &  -a_6 & a_1
 &  a_0 &  a_3 &  -a_2  \\   a_6  & -a_7 &  -a_4 &  a_5 &  a_2 &  -a_3 &  a_0 &  a_1   \\ a_7  & a_6 &  -a_5 &  -a_4 &  a_3 &  a_2 &  -a_1 &  a_0   \end{array} \right],   \eqno(2.2)
$$

\medskip

\noindent {\bf  Theorem 2.1.} \, {\em  Let $x = x_0 + x_1e_1 + \cdots +
x_7e_7 \in \Bbb O$, and denote $ \overrightarrow{x} = [ x_0, \ x_1, \ \cdots , \
x_7]^T,$  called the vector representation of $ x $. Then
$$
\overrightarrow{ax} = \omega (a) \overrightarrow{x} \eqno(2.3)
$$ 
holds for $ a, \ x \in \Bbb O$. }

\medskip

\noindent {\bf Proof.} \, Write  $ a, \, x \in \Bbb O$ as $ a = a' +  a''e,
\ x = x' +  x''e,$ where $ a',  \  a'',  \  x', \   x'' \in \Bbb H$.
We know by Eq.(1.3) that $ ax =  ( a'x'  -   \overline{x''}a'' ) +
(  x''a' +  a''\overline{x'})e.$ Thus it follows by Eq.(1.20) that
\begin{eqnarray*} 
  \overrightarrow{ax}  =  \left[ \begin{array}{c}
  \overrightarrow{ a'x'  -   \overline{x'' }a'' } \\
  \overrightarrow{  x''a' +  a''\overline{x'} } \end{array} \right]
& = & \left[ \begin{array}{c}  \overrightarrow{ a'x'}  -
\overrightarrow{ \overline{x'' }a'' } \\
 \overrightarrow{  x''a' } +  \overrightarrow{  a''\overline{x'} }  \end{array} \right] \\
& = &  \left[ \begin{array}{c}  \phi(a') \overrightarrow{x' }  -  \tau (a'') K_4\overrightarrow{x''} \\ 
\tau (a') \overrightarrow{ x'' } +  \phi(a'')K_4\overrightarrow{x'}  \end{array} \right] \\
& = & \left[ \begin{array}{cc}  \phi(a')  &  -  \tau (a'') K_4 \\  \phi(a'')K_4 & \tau (a')   \end{array} \right]  \left[ \begin{array}{c}   \overrightarrow{x'}  \\ \overrightarrow{ x'' } \end{array} \right], 
\end{eqnarray*} 
as required for Eq.(2.3).  \qquad  $ \Box $   

\medskip

\noindent {\bf  Theorem 2.2.} \, {\em  Let $a \in  \Bbb O$ be given.
Then
$$ 
aE_8 = E_8 \omega (a), \ \ and  \ \  E_8^*a  = \omega (a)E_8^*,  \eqno (2.4) 
$$
where $ E_8 := [\, 1, \  e_1, \ \cdots, \  e_7 \, ],$ and $
E_8^* := [\, 1, \ -e_1, \ \cdots, \  -e_7\, ]^T.$ }

\medskip

\noindent {\bf Proof.} \, Follows from a direct verification.
\qquad  $ \Box $  

\medskip

We can also introduce from Eq.(2.1) another matrix representation for an
octonion as follows.  

\medskip

\noindent {\bf Definition 2.2. }  \, Let $ a = a' +  a''e = a_0 + a_1e_1 +
\cdots + a_7e_7 \in \Bbb O$ be given, where $ a', \ a'' \in \Bbb H$.  Then we call the $ 8 \times 8 $ real
matrix
$$ 
\nu(a):= K_8 \omega^T (a) K_8 = \left[ \begin{array}{cc} \tau(a') &  - \phi( \overline{ a''}) 
   \\   \phi(a'')  & \tau( \overline{a' }) \end{array} \right],  \eqno (2.5)
$$     
the right matrix representation of $ a$,  where $ K_8 =
{\rm diag}( \, K_4, \ I_4  \,), $ an orthogonal matrix. Written in an
explicit form,
$$ 
\nu (a) =  \left[ \begin{array}{crrrrrrr}  a_0  & -a_1 &  -a_2 &  -a_3 &  -a_4 &  -a_5 &  -a_6 &  -a_7 
 \\   a_1 & a_0 &  a_3 &  -a_2 & a_5 &  -a_4 &  -a_7 & a_6  \\ a_2 &  -a_3 & a_0 &  a_1 &  a_6 &  a_7 &  -a_4 & -a_5  \\ a_3 &  a_2 &  -a_1 & a_0  &  a_7 &  -a_6 & a_5 & -a_4  \\   a_4  & -a_5 &  -a_6 &  -a_7 &  a_0 &  a_1 &  a_2 &  a_3  \\   a_5  & a_4 &  -a_7 &  a_6 & -a_1 &  a_0 &  -a_3 &  a_2  \\   a_6  & a_7 &  a_4 &  -a_5 &  -a_2 &  a_3 &  a_0 &  -a_1   \\ a_7  & -a_6 &  a_5 &  a_4 &  -a_3 &  -a_2 &  a_1 &  a_0  \end{array} \right].   \eqno(2.6)
$$

\noindent {\bf Theorem 2.3.} \, {\em  Let $ a, \, x \in \Bbb O$ be given. 
Then 
$$
\overrightarrow{xa} = \nu(a) \overrightarrow{x} \eqno(2.7)
$$ 
holds. } 

\medskip

\noindent {\bf Proof.} \ Write  $ a, \, x \in \Bbb O$ as $ a = a' +  a''e, \,
 x = x' +  x''e,$ where $ a',  \,  a'',  \,  x', \,   x'' \in \Bbb H$. 
 we know by (1.3)  that
 $ xa =  ( \, x'a'  -   \overline{a''}x'' \,) +
 (  \, a''x' +  x''\overline{a'} \, )e.$ Thus we find by Eq.(1.20) that
\begin{eqnarray*} 
\overrightarrow{xa} = \left[ \begin{array}{c} \overrightarrow{ x'a' - \overline{a''}x''  } \\ \overrightarrow{  a''x' +  x''\overline{a'} } \end{array} \right] & = & \left[ \begin{array}{c}  \overrightarrow{x'a'} - \overrightarrow{ \overline{a'' }x'' } \\ \overrightarrow{ a''x'} +  \overrightarrow{  x''\overline{a'} } \end{array} \right]= \left[ \begin{array}{c}  \tau(a')\overrightarrow{ x'}  -  \phi( \overline{ a''}) 
\overrightarrow{ x''} \\  \phi(a'') \overrightarrow{ x'}  + \tau( \overline{a'})\overrightarrow{ x''}  \end{array} \right] = \left[ \begin{array}{cc} \tau(a') &  - \phi( \overline{ a''}) \\   \phi(a'')  & \tau( \overline{a'}) \end{array} \right] \left[ \begin{array}{c}  \overrightarrow{ x'}  \\ \overrightarrow{ x''}  \end{array} \right], 
\end{eqnarray*} 
as required for Eq.(2.7).  \qquad   $ \Box $  

\medskip

\noindent {\bf Theorem 2.4.} \, {\em  Let $a \in \Bbb O$ be given. Then 
$$ 
aF_8 = F_8 \nu^T(a), \ \ and  \ \  F_8^*a  = \nu^T (a)F_8^*,  \eqno (2.8) 
$$
where $ F_8 := [ \, 1, \  -e_1, \ \cdots, \  -e_7 \, ]$ and
$ F_8^* := [ \, 1, \ e_1, \ \cdots, \  e_7 \, ]^T.$ } 

\medskip

\noindent {\bf Proof.} \, Follows from a direct verification.
 \qquad  $ \Box $  

\medskip

Observe from Eqs.(2.1) and (2.5) that the two real matrix representations of
an octonion $ a =  a' + a''e $ are in fact constructed by the real matrix
representations of two quaternions $ a'$ and  $a''$. Hence the operation
properties for the two matrix representations of octonions  can easily be
established through the results in Lemmas 1.2 and 1.3.
 
\medskip

\noindent {\bf  Theorem 2.5.} \, {\em  Let $a, \ b \in \Bbb O,\ \lambda
\in  \Bbb R$  be given. Then
   
{\rm (a)} \ $ a=b \Longleftrightarrow \omega(a) = \omega(b).$ 

{\rm (b)} \ $ \omega(a+b)= \omega(a) + \omega(b),
 \qquad \omega( \lambda a)=  \lambda\omega(a), \qquad  \omega(1) = I_8.$

{\rm (c)} \ $  \omega(\overline{a})= \omega^T(a). $  }

\medskip

\noindent {\bf Proof.} \, Follows from a direct verification. \qquad $\Box$ 

\medskip

\noindent {\bf  Theorem 2.6.}  \, {\em  Let $a, \, b \in \Bbb O,\, \lambda
\in \Bbb R$ be given. Then

{\rm (a)} \ $ a=b  \Longleftrightarrow \nu(a) = \nu(b).$

{\rm (b)} \ $ \nu(a+b)= \nu(a) + \nu(b), \qquad
 \nu( \lambda a)=  \lambda\nu(a), \qquad  \nu(1) = I_8.$

{\rm (c)} \ $  \nu (\overline{a})= \nu^T(a). $  }  

\medskip

\noindent {\bf Proof.} \ Follows from a direct verification. \qquad $\Box$ 

\medskip

\noindent {\bf  Theorem 2.7.}  \, {\em  Let $a \in \Bbb O$ be given.
Then
$$ 
a =  \frac{1}{8} E_8 \omega(a)E_8^*, \ \ and \ \
a =  \frac{1}{8} F_8 \nu^T(a)F_8^*, \eqno( 2.9)
$$ 
where $ E_8, \ E_8^*, \  F_8$ and  $F_8^* $ are as in Eqs.{\rm (2.4)} and
{\rm (2.8)}. }

\medskip

\noindent {\bf Proof.} \, Note that $\omega (a)$ and $\nu(a)$ are
real matrices. Thus we  get from Eqs.(2.4) and (2.8) that
$$ 
E_8(E_8^*a ) = E_8 [ \omega(a)E_8^* ] = E_8 \omega(a)E_8^*, \ \ and \ \
 F(F_8^*a ) = F_8 [ \nu^T(a)F_8^* ] = F_8 \nu^T(a)F_8^*.
$$ 
On the other hand, note that $\Bbb O$ is alternative. It follows that
$$ 
E_8(E_8^*a ) = a - e_1(e_1a) - \cdots - e_7(e_7a)  =  a - e_1^2a  - \cdots - 
 e_7^2a = 8a, 
$$ 
and
$$
F_8(F_8^*a ) = a - e_1(e_1a) - \cdots -  e_7(e_7a)  =  a - e_1^2a  - \cdots - 
 e_7^2a = 8a.
$$ 
Thus we have Eq.(2.9).  \qquad  $ \Box $   

\medskip

\noindent {\bf  Theorem 2.8.} \, {\em  Let $a  \in \Bbb O$ be given.
 Then
 $$
 {\rm det\,}[\omega(a)]= {\rm det\,} [ \nu (a)] =  |a|^8. \eqno (2.10)  
$$   }
\noindent {\bf Proof.} \ Write $ a = a' +  a''e $. Then  we easily find 
 by Eqs.(1.21) and (2.5) that
\begin{eqnarray*} 
{\rm det\,}[\omega(a)]= {\rm det\,} [ \nu (a)]   =
\left| \begin{array}{cc} \tau(a') &  - \phi( \overline{ a''})
   \\   \phi(a'')  & \tau( \overline{a'}) \end{array} \right|
   & = & {\rm det\,}[ \, \tau(a')  \tau( \overline{a' }) +
   \phi(a'')\phi( \overline{ a''}) \,] \\
&= & {\rm det\,} [ \, \tau( \overline{a'} a') + \phi(a''\overline{ a''}) \,] \\
&= & {\rm det\,} [ \, |a'|^2I_4 +  |a''|^2I_4 \, ] \\
& = & ( \, |a'|^2 +  |a''|^2 \,)^4 = |a|^8,
\end{eqnarray*} 
as required for Eq.(2.10).  \qquad  $ \Box $    

\medskip

\noindent {\bf Theorem 2.9.} \, {\em  Let $a  \in \Bbb O$ be given.
Then the two matrix representations of $ a $ satisfy the following three
identities
$$
 \omega(a^2)  =  \omega^2(a), \qquad    \nu (a^2) = \nu^2(a), \qquad
 \omega(a)\nu(a) = \nu(a) \omega(a). \eqno (2.11)  
$$   } 

\noindent {\bf Proof.} \, Applying Eqs.(2.3) and (2.7) to the both sides of
the three identities in Eq.(1.9) leads to
$$ 
\omega^2(a)\overrightarrow{b} =  \omega(a^2)\overrightarrow{b},  \ \ \
\nu^2(a)\overrightarrow{b} =  \nu(a^2)\overrightarrow{b}, \ \ \
  \omega(a)\nu(a)\overrightarrow{b} = \nu(a) \omega(a)\overrightarrow{b}.
$$ 
Note that $ \overrightarrow{b}$ is an arbitrary $ 8 \times 1 $ real vector
when $ b $ runs over  $ \Bbb O$. Thus Eq.(2.11) follows.  \qquad  $ \Box $ 

\medskip

\noindent {\bf Theorem 2.10.} \, {\em  Let $a  \in \Bbb O$ be given
 with $ a \neq 0 $.  Then
$$
 \omega(a^{-1})= \omega^{-1}(a), \ \ \ and  \ \ \  \nu(a^{-1})= \nu^{-1}(a).
    \eqno (2.12)
$$ } 
\noindent {\bf Proof.} \, Note from Eqs.(1.10) and (1.11) that  
$$ 
a^{-1} = \frac{ \overline{a}}{|a|^2} = \frac {1}{ |a|^2}
[\, 2({\rm Re}\,a) - a \, ]
$$ 
and
$$
a^2 - 2{\rm Re\,} a + |a|^2 = 0.
$$ 
Applying Theorems 2.5 and 2.6, as well as the first two equalities in
Eq.(2.11) to the both sides of the above two equalities,
 we obtain
$$
\omega(a^{-1}) = \frac {1}{ |a|^2} [\, 2({\rm Re}\,a)I_8  - \omega(a) \,], 
 \qquad
\nu(a^{-1}) = \frac {1}{ |a|^2} [\, 2({\rm Re}\,a)I_8  - \nu(a) \,]
$$ 
and
$$
\omega^2(a) - 2({\rm Re}\,a )\omega(a) + |a|^2 I_8 = 0, \qquad
 \nu^2(a) - 2({\rm Re}\,a )\nu(a) + |a|^2 I_8 = 0.
$$ 
Contrasting them yields Eq.(2.12). \qquad  $ \Box$

\medskip

Because $\Bbb O$ is non-associative, the operation properties
$ \omega(ab) =  \omega(a) \omega(b)$ and $ \nu(ab) = \nu(b) \nu(a)$
 do not hold in general, otherwise $\Bbb O$ will be
algebraically isomorphic to or algebraically anti-isomorphic to an
associative matrix algebra over $ \Bbb R$, this is impossible. Nevertheless, 
some other kinds of identities on the two real matrix representations of 
octonions  can  still be established from the  
identities in Lemma 1.4(a)---(d). 

\medskip

\noindent {\bf Theorem 2.11.} \, {\em Let $a, \, b  \in \Bbb O$ be given.
 Then  their matrix representations satisfy the  following two identities 
$$ 
\omega(aba) = \omega(a) \omega(b) \omega(a),  \  \ and \ \
 \nu(aba) =\nu(a)\nu(b) \nu(a).  \eqno (2.13)
$$ } 
\noindent {\bf Proof.} \, Follows from applying Eqs.(2.3) and (2.7) to the
Moufang identities in Lemma 1.4(b) .  \qquad $ \Box $  

\medskip

\noindent {\bf Theorem 2.12.} \, {\em Let $a, \, b  \in \Bbb O$ be given.
  Then  their matrix representations satisfy  the  following identities
$$
 \displaylines{
\hspace*{2cm} \omega (ab) + \omega (ba)   = \omega (a) \omega (b) + \omega (b) \omega (a),
 \hfill (2.14)
 \cr
 \hspace*{2cm} \nu (ab) + \nu (ba) = \nu (a)\nu (b) + \nu (b)\nu (a),  \hfill (2.15)
 \cr
 \hspace*{2cm}
 \omega (ab) + \nu (ab)  = \omega (a) \omega (b) + \nu (b)\nu (a),  \hfill (2.16)
 \cr
 \hspace*{2cm}
 \omega (a) \nu (b)  +  \omega (b) \nu (a) =   \nu (a) \omega (b) + \nu (b)
 \omega (a),  \hfill (2.17)
 \cr
 \hspace*{2cm}
 \omega (ab) = \omega (a) \omega (b) +  \omega (a) \nu (b) - \nu (b) \omega (a),
 \hfill (2.18)
 \cr
 \hspace*{2cm}
 \nu (ab)  = \nu (b) \nu (a) +   \omega(b) \nu (a) - \nu (a) \omega (b). 
 \hfill (2.19)
 \cr }
$$ }
\noindent {\bf Proof.} \, The identities in Lemma 1.4(d) can  clearly  be
 written as the following six identities 
$$ 
(ab)x - a(bx) = -(ba)x + b(ax),  \qquad (xa)b - x(ab) = -(xb)a + x(ba), 
$$ 
$$ 
(ab)x - a(bx) = -(bx)a + b(xa), \qquad  (ab)x - a(bx) = -(xa)b + x(ab),  
$$ 
$$ 
(ab)x - a(bx) =  -(ax)b + a(xb), \qquad (xa)b - x(ab) = -(ax)b + a(xb).
$$ 
Applying Eqs.(2.3) and (2.7) to the both sides of the above 
identities,  we obtain
$$
\displaylines{
\hspace*{2cm} 
 [\, \omega(ab) - \omega(a) \omega(b) \, ]\overrightarrow{x} = [ \, - \omega(ba) + \omega(b) \omega(a) \, ] 
\overrightarrow{x},  \hfill 
\cr
\hspace*{2cm}
  [ \, \nu(b)\nu(a) - \nu(ab) \, ]\overrightarrow{x} = [ \, -\nu(a)\nu(b)  + \nu(ba) \, ] \overrightarrow{x},  \hfill 
\cr
\hspace*{2cm}
  [ \, \nu(b)\omega(a) - \omega(a)\nu(b)  \, ]\overrightarrow{x} =  [\,  -\nu(a)\omega(b)  + \omega(b) \nu(a) \, ] \overrightarrow{x},  \hfill 
\cr
\hspace*{2cm} 
 [ \, \omega(ab) - \omega(a) \omega(b)  \, ]\overrightarrow{x} =  [\,  -\nu(b)\nu(a)  +  \nu(ab)  \, ]\overrightarrow{x}, \hfill 
\cr
\hspace*{2cm} 
  [ \, \omega(ab) - \omega(a) \omega(b) \, ]\overrightarrow{x} = [\, -\nu(b)\omega(a)  + \omega(a) \nu(b)\, ] \overrightarrow{x}, \hfill 
 \cr
\hspace*{2cm}
 [ \, \nu(b)\nu(a) - \nu(ab) \, ]\overrightarrow{x} = [\, -\nu(b)\omega(a)  + \omega(a)\nu(b)\, ] \overrightarrow{x}.  \hfill 
\cr }
$$
Notice  that $\overrightarrow{x}$ is an arbitrary real $ 8 \times1 $
 real matrix when $ x $ runs over $\Bbb O$. Therefore Eqs.(2.14)---(2.19)
  follow.   \qquad  $ \Box $  

\medskip

\noindent {\bf Theorem 2.13.} \, {\em Let $a, \, b  \in \Bbb O$ be given
with $ a \neq 0, \,  b \neq 0.$   Then  their matrix representations satisfy
 the  following two identities
$$ 
\omega(ab)  = \nu(a)[\, \omega(a) \omega(b) \,] \nu^{-1}(a),  \ \  and \ \
\nu(ab)  = \omega(b)[\, \nu(b) \nu(a)\,] \omega^{-1}(b). \eqno (2.20)
$$
which imply that 
$$
 \omega(ab) \sim \omega(a)\omega(b), \ \  and \ \ \nu(ab) \sim \nu(b)\nu(a). 
\eqno (2.21)
$$}

\noindent {\bf Proof.} \ Applying Eqs.(2.3) and (2.7) to  the both sides of
the two  identities in Lemma 1.4(c),
we obtain
$$ 
\omega(ab) \nu(a) \overrightarrow{x} = \nu(a)\omega(a) \omega(b)
\overrightarrow{x}, \ \ and \ \   \nu(ab) \omega(b) \overrightarrow{x} =
\omega(b)\nu(b) \nu(a) \overrightarrow{x},
$$
which are obviously equvalent to Eq.(2.20).  \qquad $ \Box $  

\medskip





Note from Eqs.(2.3) and (2.7) that any linear equation of the form
$ ax - xb =c$ over $\Bbb O$ can equivqlently be written as
$ [\, \omega (a)- \nu (b) \,] \overrightarrow{x} = \overrightarrow{a}$,
which is a linear equation over $\Bbb R$.  Thus it is necessary to
consider the operation properties of the matrix $ \omega(a ) - \nu( b)$, especially the determinant 
of $ \omega(a ) - \nu( b)$ for any $ a, \ b \in \Bbb O$. Here we only list the expression  of the  
determinant of $ \omega(a ) - \nu( b)$. Its proof is quite tedious and is, therefore, omitted
 here.

\medskip 

\noindent {\bf Theorem 2.14.} \, {\em Let $ a, \, b \in \Bbb O$ be given and define
$ \delta( a, \ b ) :=  \omega(a )- \nu( b ).$ Then 
$$
\displaylines{
\hspace*{0.2cm}
{\rm det\,}[ \delta( a, \ b ) ] = | a - \overline{b}|^4 [ \, s^2 + ( \, |{\rm
Im}\,a| -|{\rm Im}\,b| \, )^2 \, ]
 [ \, s^2 + ( \ |{\rm Im}\,a| + |{\rm Im}\,b| \, )^2 \,] \hfill (2.22)
\cr
\hspace*{0.2cm} 
{\rm det\,}[ \delta( a, \ b ) ] = 
( \, s^2 + |{\rm Im}\,a + {\rm Im}\,b|^2  \,)^2[ s^4 +  2s^2 (|{\rm Im}\,a|^2 +|{\rm Im}\,b|^2 ) +
( \, |{\rm Im}\,a|^2 -|{\rm Im}\,b|^2 \, )^2 ],  \hfill (2.23)
  \cr  }
$$
where $ s = {\rm Re\,}a -{\rm Re\,}b.$ The characteristic polynomial of 
$\delta( a, \ b )$ is
$$
\displaylines{
\hspace*{0cm}
| \: \lambda I_8 -  \delta( a, \ b ) \: | \hfill
\cr
\hspace*{0cm}
= [ \: ( \lambda - s \: )^2 +|{\rm Im}\,a + {\rm Im}\,b|^2 ]^2
[ \: ( \lambda - s \: )^2 + ( |{\rm Im}\,a| -|{\rm Im}\,b| \ )^2 \: ][\: ( \lambda - s \: )^2 + ( \: |{\rm Im}\,a| + |{\rm Im}\,b| \: )^2 \: ]. \hfill (2.24)
\cr }
$$
In particular$,$ if  $ {\rm Re\,}a ={\rm Re\,}b$ and $  |{\rm Im}\,a|
= |{\rm Im}\,b|,$ but $ a \neq \overline{b},$ then
$$ 
 {\rm rank \,}\delta( a, \, b ) = 6.  \eqno (2.25) 
$$ }
\noindent {\bf Theorem 2.15.} \, {\em Let $a, \, b  \in \Bbb O$ be given
with $ a \neq 0$ and $ b \neq 0.$   Then $ \delta( a, \ b ) =  \omega(a ) -
\nu( b )$ is a real normal matrix over $\Bbb R,$ that is$,$
$ \delta( a, \ b ) \delta^T( a, \ b ) = \delta^T( a, \ b ) 
\delta( a, \ b )$. }

\medskip

\noindent {\bf Proof.} \, Follows from
\begin{eqnarray*} 
  \delta( a,  \ b ) + \delta^T( a, \ b ) & = & \omega( a ) - \nu(b) +
   \omega^T( a ) - \nu^T(b) \\
& = &  \omega( a ) - \nu(b) + \omega( \overline{ a } ) - \nu( \overline{ b }) \\
& = &  \omega( a + \overline{ a } ) - \nu(  b + \overline{ b }) =
2( {\rm Re}\,a - {\rm Re}\,b )I_8. \qquad \Box 
\end{eqnarray*} 

\noindent {\bf Theorem 2.16.} \, {\em Let $a \in \Bbb O$ be given with
$ a \notin \Bbb R$. Then
 $$ 
\delta^3( a, \ a ) = -4|{\rm Im}\,a|^2 \delta( a, \ a ),  \eqno (2.26) 
$$ 
and $  \delta( a, \ a)$ has a generalized inverse as follows
$$ 
 \delta^{-}( a, \ a ) = - \frac{1}{4|{\rm Im}\,a|^2} \delta( a, \ a ).
 \eqno (2.27)
$$ } 
\noindent { \bf Proof.} \, Observe that $ \delta( a, \ a) = \omega( a ) -
\nu(a) =  \omega({\rm Im}\,a ) - \nu({\rm Im}\,a)$ and
$ ({\rm Im}\,a )^2 = - |{\rm Im}\,a|^2$. Thus we find that
\begin{eqnarray*} 
\delta^2( a, \ a ) & = & [\, \omega({\rm Im}\,a ) - \nu({\rm Im}\,a) \,]^2 \\
  & = & [\,  \omega^2({\rm Im}\,a ) - 2 \omega({\rm Im}\,a )
  \nu({\rm Im}\,a)  +  \nu^2({\rm Im}\,a) \,]\\
  & = & [\,  \omega(({\rm Im}\,a )^2) - 2 \omega({\rm Im}\,a )
  \nu({\rm Im}\,a)  +  \nu(({\rm Im}\,a)^2) \,]\\
 & = & -2 [\,  |{\rm Im}\,a|^2I_8 + \omega({\rm Im}\,a ) \nu({\rm Im}\,a) \,],
 \end{eqnarray*}
and  
\begin{eqnarray*} 
\delta^3( a, \ a ) & = & -2 [ \, |{\rm Im}\,a|^2I_8 + \omega({\rm Im}\,a )
\nu({\rm Im}\,a) \,][\, \omega( {\rm Im}\,a ) - \nu( {\rm Im}\,a) \,]  \\
 & = &   - 4|{\rm Im}\,a|^2 [\, \omega({\rm Im}\,a ) -
 \nu({\rm Im}\,a) \,] = -4|{\rm Im,}a|^2 \delta(a,\ a),
 \end{eqnarray*}
as required  for Eq.(2.26).   \qquad $ \Box $ \\

\noindent  {\Large {\bf 3. \  Some linear equations over $ \Bbb O$}} \\

\noindent The matrix expressions of octonions and their properties 
 introduced in Section 2 enable us to easily deal with various 
problems related to octonions. One of the most fundamental topics on octonions 
is concerning solutions of various linear equations over $ \Bbb O$. 
In this section, we shall give a complete discussion for this problem. Our first 
result is concerning the linear equation $ ax = xb$, which was examined by the author
in \cite{Ti1}.     

\medskip

\noindent {\bf  Thoerem 3.1}\cite{Ti1}. \, {\em  Let $a = a_0 + a_1e_1 + \cdots +
a_7e_7, \,  b = b_0 + b_1e_1 + \cdots + b_7e_7  \in \Bbb O$ be given.
Then the linear equation $ ax = xb$ has a nonzero solution if and only if  
$$
{\rm Re}\,a  ={\rm Re}\,b  \ \ \  and  \ \ \ | {\rm Im}\,a | = |
{\rm Im}\,b |. \eqno (3.1)
 $$ 

{\rm (a)} \ In that case$,$ if $ b \neq \overline{a},$
i. e.$,$  $ {\rm Im}\,a  + {\rm Im}\,b \neq 0,$  then the general solution of
 $ ax = xb$ can be expressed as
$$ 
 x = ({\rm Im\,}a )p + p({\rm Im\,}b),  \eqno (3.1)
$$ 
where  $ p \in {\cal A}(a, \, b )$, the subalgebra generated by $ a $ and
$ b $, is arbitrary. or equivalently
 $$
 x = \lambda_1( \, {\rm Im\,}a  + {\rm Im\,}b \,) + \lambda_2 [\,
 |{\rm Im\,}a | \,|{\rm Im\,}b| -( {\rm Im\,}a )({\rm Im\,}b ) \,],
 \eqno (3.2)
$$ 
where $ \lambda_1, \, \lambda_2 \in \Bbb R$ are arbitrary. 

{\rm (b)}  \ If $ b = \overline{ a },$ then the general solution of
 $ ax = xb$ is
$$
 x =  x_1e_1 + x_2e_2  + \cdots +  x_7e_7, \eqno (3.3)
$$
where  $x_1$---$x_7$ satisfy $ a_1x_1 + a_2x_2  + \cdots + a_7x_7 = 0.$  }
 
\medskip

The correctness of this result can be directly  verified by substitution.

\medskip
 
Based on the equation $ ax = xb$,  we can define the similarity of two octonions.
Two octonions are said to be {\em  similar} if there is a nonzero
$ p \in \Bbb O$  such that $ a = pbp^{-1}$, which  is  written
as $ a \sim b$. Theorem 3.1 shows  that two octonions are
similar  if and only if  $  {\rm Re\,}a  = {\rm Re\,}b$  and
$ | {\rm Im\,}a | = | {\rm Im\,}b |$. Thus  the similarity defined here is
also an equivalence relation on octonions. In addition, we have the
following.

\medskip

\noindent {\bf  Theorem 3.2.} \, {\em  Let $a, \, b \in \Bbb O$ be given
with  $ b \neq \overline{a}$. Then
$$ 
 a \sim b  \Longleftrightarrow  \omega( a )  \sim \omega( b ). \eqno (3.4)  
$$ } 
\noindent {\bf Proof.} \ Suppose first that  $ a \sim b$. Then it follows by
Eq.(1.11)  that
$$ 
a^2 - 2 ({\rm Re\,}a)a = - |a|^2 = - |b|^2 = b^2 - 2 ({\rm Re\,}b)b.
$$
Applying Theorem 2.5(a) and Eq.(2.11) to the both sides of the above 
equality  and  we get
$$   
 \omega^2(a) -  2 ({\rm Re\,}a)\omega(a) = \omega^2(b) -  2 ({\rm Re\,}b)\omega(b).
$$ 
Thus  
$$   
 \omega^2(a) + \omega(a)\omega(b) -  2 ({\rm Re\,}a)\omega(a)
 = \omega^2(b) + \omega(a)\omega(b) -
 2 ({\rm Re\,}b)\omega(b), 
$$ 
which is equivalent to 
$$   
 \omega(a)[ \,  \omega(a) + \omega(b) -  2 ({\rm Re\,}a)I_8 \,] =
 [ \,  \omega(a) + \omega(b) -  2 ({\rm Re\,}b)I_8 \,] \omega(b),
$$
or simply 
 $$ 
\omega(a) \omega( {\rm Im\,}a + {\rm Im\,}b )
= \omega( {\rm Im\,}a + {\rm Im\,}b )\omega(b).
$$    
Note that $ {\rm Im\,}a + {\rm Im\,}b \neq 0$. Thus
 $ \omega( {\rm Im\,}a + {\rm Im\,}b )$ is invertible.  The above equality 
shows that $ \omega( a ) \sim \omega( b ) $.  Conversely, if $ \omega( a ) \sim
\omega( b )$, then trace\,$\omega( a )=$ trace\,$\omega( b )$ and
$ |\omega(a)| =  |\omega(b)|$, which are equivalent to Eq.(3.1).
 \qquad $ \Box $  

\medskip

Next we consider some  nonhomogeneous linear equations over $ \Bbb O$.
 
\medskip

\noindent {\bf  Theorem 3.3.} \, {\em  Let $a, \, b \in \Bbb O$ be given
 with $ a \notin \Bbb R$. Then the linear equation $ ax -xa = b$ has a 
solution in $ \Bbb O$ if and only if The equality $ ab = b\overline{a}$ holds. In this case$,$ the general solution of $ ax - xa = b$ is
$$
 x = \frac{1}{4|{\rm Im}\,a |^2}( \,ba - ab \, ) +
 p - \frac{1}{ |{\rm Im}\,a|^2}({\rm Im}\,a )p({\rm Im}\,a ),
  \eqno (3.5)
$$ 
where $ p \in \Bbb O$ is arbitrary. }

\medskip

\noindent {\bf Proof.} \,  According to Eqs.(2.3) and (2.7), the equation
$ ax - xa = b$   can  equivalently be written as
$$ 
 [  \, \omega( a ) -  \nu(a)  \,] \overrightarrow{ x } =\delta( a, \ a )
 \overrightarrow{ x }
= \overrightarrow{b }. \eqno (3.6)
$$
This equation is solvable if and only if $\delta( a, \ a )
\delta^{-}( a, \ a )\overrightarrow{b} = \overrightarrow{b}.$
In that case, the general  solution of Eq.(3.6) can be expressed as
$$
 \overrightarrow{ x } = \delta^{-}( a, \ a ) \overrightarrow{c} +
 2[ \,I_8 -   \delta^{-}( a, \, a ) \delta( a, \, a ) \, ]\overrightarrow{ p },
$$
where $ \overrightarrow{ p }$ is an arbitrary real vector. Substituting
$$ 
 \delta^{-}( a, \ a ) = - \frac{1}{4|{\rm Im}\,a|^2}\delta(a, \ a) , \ \ and
 \ \  \delta^2( a, \ a )= -2 [\, |{\rm Im}\,a|^2 + \omega({\rm Im}\,a )
 \nu({\rm Im}\,a) \,]
$$
in the above two equalities and then returning them to octonion forms by
Eqs.(2.3) and (2.7) produce the equality in Part (b) and
Eq.(3.5).  \qquad  $ \Box $

\medskip

\noindent {\bf Theorem 3.4.}  \, {\em  Let $ a = a_0 + a_1e_1 + \cdots +
a_7e_7 , \, b = b_0 + b_1e_1 + \cdots + b_7e_7  \in \Bbb O$ be given
with $ a \notin \Bbb R.$  Then the equation
$$ 
ax -x\overline{a} = b  \eqno (3.7)
$$
has a solution if and only if there exist $  \lambda_0, \, \lambda_1 \in
 \Bbb R $ such that
$$ 
 b = \lambda_0 + \lambda_1a , \eqno (3.8)
$$ 
in which case$,$ the general solution of Eq.{\rm (3.7)} is 
$$
 x =\frac{\lambda_1}{2} + x_1e_1 + \cdots + x_7e_7, \eqno (3.9)
$$ 
where $ x_1 $---$x_7$ satisfy 
$$ 
a_1x_1 + \cdots +  a_7x_7 = - \frac{1}{2}{\rm Re\,}b. \eqno (3.10)
$$ }
\noindent {\bf Proof.} \, According to  Eqs.(2.3) and (2.7), the equation
(3.7) is equivalent to
$$ 
[\, \omega( a ) -  \nu(\overline{a}) \,] \overrightarrow{ x } =
\delta( a, \ \overline{a} ) \overrightarrow{ x } = \overrightarrow{b },
\eqno (3.11)
$$ 
namely 
$$ 
  \left[ \begin{array}{cccc} 0 & -2a_1  &  \cdots &  -2a_7  \\ 2a_1  & 0
  &  \cdots &  0 \\ \vdots &
\vdots & \ddots & \vdots  \\ 2a_7 & 0 &  \cdots &  0  \end{array} \right]
\left[ \begin{array}{c} x_0  \\
 x_1 \\ \vdots \\ x_7  \end{array} \right]
 = \left[ \begin{array}{c} b_0  \\ b_1 \\ \vdots \\ b_7  \end{array}
 \right].
$$ 
Obviously, this equation is solvable if and only if there is
a $ \lambda_1 \in  \Bbb R$ such that
$$ 
b_1 = \lambda_1a_1, \ \ \  b_2 = \lambda_1a_2, \ \ \  \cdots, \ \ \
b_7 = \lambda_1a_7,
$$ 
i. e., ${\rm Im\,b} = \lambda_1{\rm Im\,a}$, which is equivalent to
Eq.(3.8). In that case, the solution to  $ x_0 $ is $ x_0
= \frac{\lambda_1}{2}$, and $ x_1$---$x_7$ are determined by Eq.(3.9).
 \qquad $ \Box $  

\medskip

Next we consider the linear equation 
$$ 
ax - xb = c \eqno (3.12)
$$
under the condition $ a \sim b $. Clearly Eq.(3.12) is equivalent to 
$$ 
 [\, \omega( a ) -  \nu(b) \,] \overrightarrow{ x } = \delta( a, \ b )
 \overrightarrow{ x }
= \overrightarrow{c}. \eqno (3.13)
$$
Under  $ a \sim b $, we know by Theorem 3.3 that $ ax = xb $ has a nonzero
 solution. Hence $ \delta(a, \ b )$  is singular  under  $ a \sim b $. In
 that case, Eq.(3.12) is solvable if and only if
 $$ 
 \delta( a, \ b ) \delta^{-}( a, \ b )\overrightarrow{c} =
 \overrightarrow{c}, \eqno(3.14)
$$ 
and the general solution of Eq.(3.13) is 
$$
\overrightarrow{x} = \delta^{-}( a, \ b) \overrightarrow{c} +
2[ \,I_8 -   \delta^{-}( a, \, b )
\delta( a, \, b) \, ]\overrightarrow{p}, \eqno(3.15)
$$
where $ \overrightarrow{p}$ is an arbitrary real vector. If 
 $ a$ is not not similar to  $b$. Clearly  Eq.(3.13)  has a unique solution    
$$
\overrightarrow{x} = \delta^{-1}( a, \ b) \overrightarrow{c} 
\eqno(3.16)
$$ 
Eqs.(3.15) and (3.16) show that the solvability and solution of the 
octonion equation (3.12) can be completely determined by its real adjoint 
 linear system of equations (3.13). Through the characteristic polynomial 
(2.24), one can also retern Eqs.(3.15) and (3.16) to octonion forms. But 
their expressions are quite tedious in form, and are omitted here.

\medskip

Another instinctive linear equation over $ \Bbb O$ is 
$$
a(xb) - (ax)b = c,  \eqno (3.17)
$$  
which is also equivalent to 
$$ 
(ab)x - a(bx) = c,  \eqno (3.18)
$$  
as well as 
$$ 
x(ab) - (xa)b = c,  \eqno (3.19)
$$  
because $(ab)x - a(bx) = (ab)x - a(bx) =x(ab) - (xa)b $ hold for all
$  a, \ b, \ x  \in \Bbb O$. Now applying Eqs.(2.3) and (2.7) to the both
sides of Eq.(3.17), we obtain an equivalent equation
$$ 
  [ \, \omega(a)\nu(b) - \nu(b)\omega(a) \, ] \overrightarrow{ x }  =
  \overrightarrow{ c }.
\eqno (3.20)
$$
Here we set $ \mu( a, \ b ) = \omega(a)\nu(b) - \nu(b)\omega(a) $. Then it
is easy to see that Eq.(3.20) is solvable if and only if 
$$ 
\mu( a, \ b )\mu^{-}( a, \, b )  \overrightarrow{ c } =
\overrightarrow{ c },
$$ 
where $ \mu^{-}( a, \, b )$ is a generalized inverse of $\mu( a, \, b )$.
In that case, the general solution of Eq.(3.20) is
$$ 
\overrightarrow{x} = \mu^{-}( a, \, b) \overrightarrow{c} + [ \,I_8 -
\mu^{-}( a, \, b ) \mu( a, \, b) \, ]\overrightarrow{p}, \eqno(3.21)
$$
where $ \overrightarrow{p}$ is an arbitrary real vector. Numerical computation for 
Eq.(3.21) can reveal some interesting facts on Eq.(3.17). The reader can try to find them. 


\medskip

Theoretically speaking, any kind of two-sided linear equations or systems of linear
equations  over $  \Bbb O$ can be equivalently transformed into  systems of linear equations
 over $  \Bbb R$ by the two equalities in Eqs.(2.3) and (2.7). Thus the problems related to linear 
equations over $ \Bbb O $ now have a complete resolution. \\

\noindent  {\Large {\bf 4. \ Real adjoint matrices of octonion matrices}} \\

\noindent In this section, we consider how to extend the work in Sections 2 and 
 3  to octonion matrices and use them to deal with various  octonion matrix problems. 
  Since octonion algebra is non-associative, the matrix operations
in $\Bbb O$ is much different from what we are familiar with in an
associative algebra. Even the simplest matrix multiplication rule
$ A^2A = AA^2 $ does not hold over $\Bbb O$, that is to say, 
multiplication of matrices over $\Bbb O$ is completely not associative.
 Thus nearly all the known results and methods on matrices over
 associative algebras  can hardly be  extended to matrices over
 $\Bbb O$. In that case, a unique method available to deal with matrices
 over $\Bbb O$ is to establish real matrix representations of octonion matrices, 
and then to transform matrix problems over $\Bbb O$ to various equivalent
 real matrix problems. 

\medskip

Based on the two matrix representations of octonions shown in Eqs.(2.2)
and (2.6), we now introduce two adjoints for a octonion matrix as follows.

 \medskip
   
\noindent {\bf Definition 4.1.} \, Let $ A = (a_{st}) \in \Bbb O^{ m
\times n}$ be given . Then the {\em left adjoint matrix} of $ A $ is defined
to be 
$$ 
\omega( A ) =  [ \omega (a_{st})] = \left[ \begin{array}{ccc}
\omega(a_{11}) & \cdots &
\omega(a_{1n})   \\ \vdots & &   \vdots  \\  \omega(a_{m1}) & \cdots &
\omega(a_{mn})  \end{array} \right]
\in  \Bbb R^{ 8m \times 8n}, \eqno (4.1) 
$$ 
the {\em right adjoint matrix} of $ A $ is defined to be
$$
\nu( A ) =  [ \nu(a_{ts})] = \left[ \begin{array}{ccc}  \nu(a_{11})
& \cdots &
\nu(a_{m1})   \\ \vdots & &   \vdots  \\  \nu(a_{1n}) & \cdots & \nu(a_{mn})  \end{array}\right]
\in  \Bbb R^{ 8n \times 8m}, \eqno (4.2) 
$$ 
and the {\em adjoint vector} of $ A $ is defined to be
$$
{\rm vec}A := [\, \overrightarrow{a_{11}}^T, \ \cdots, \
\overrightarrow{a_{m1}}^T,\ \overrightarrow{a_{12}}^T, \ \cdots, \
\overrightarrow{a_{m2}}^T, \ \cdots,  \  \overrightarrow{a_{1n}}^T, \
\cdots, \ \overrightarrow{a_{mn}}^T \, ]^T.   \eqno (4.3)
$$ 
 \noindent {\bf Definition 4.2.} \, Let $ A = ( A_{st} )_{ m\times n}$ and
 $ B  = ( B_{st} )_{ p\times q}$ are two block matrices over $ \Bbb R$,
 where $  A_{st}, \ B_{st} \in \Bbb R^{ 8 \times 8}$. Then the
{\em left and right block Kronecker products} of $ A $ and $ B $, denoted
respectively by $ A \widehat{\otimes}B $ and $ A \tilde{\otimes}B$, are
defined to be
$$
A\widehat{\otimes}B = \left[ \begin{array}{ccc} A_{11}\odot_L B   & \cdots &
A_{1n}\odot_L B   \\
 \vdots & \ddots & \vdots  \\ A_{m1}\odot_L B  & \cdots & A_{mn}\odot_L B
  \end{array} \right] \in  \Bbb R^{8mp \times 8nq }
, \eqno (4.4)
 $$ 
and 
$$
A \tilde{\otimes}B = \left[ \begin{array}{ccc} A\odot_R B_{11}
& \cdots & A\odot_R B_{1q}    \\ \vdots & \ddots & \vdots  \\
A\odot_R B_{p1} & \cdots &   A\odot_R B_{pq}   \end{array} \right]
\in  \Bbb R^{8mp \times 8nq }, \eqno (4.5)
 $$ 
where 
$$
A_{st} \odot_L B = \left[ \begin{array}{ccc} A_{st} B_{11}   & \cdots &
A_{st} B_{1q}   \\ \vdots & \ddots & \vdots  \\ A_{st} B_{p1}  & \cdots &
A_{st} B_{pq}   \end{array} \right]
 \in  \Bbb R^{8p \times 8q},  \eqno (4.6)
 $$ 
$$
A \odot_R B_{st} = \left[ \begin{array}{ccc} A_{11} B_{st}   & \cdots &  A_{1n} B_{st}   \\
 \vdots & \ddots & \vdots  \\ A_{m1} B_{st}  & \cdots & A_{mn} B_{st}   \end{array} \right]
 \in  \Bbb R^{8m \times 8n}.   \eqno (4.7)
 $$ 

Noticing the equality (2.5), we see the two adjoint matrices $ \omega( A )$
and  $ \nu( A) $  of an octonion matrix $ A $ satisfy the following equality
$$ 
 \nu( A) = K_{8n}\omega^T( A ) K_{8m}, \eqno (4.8)
$$ 
where 
$$ 
K_{8t} = {\rm diag}( \, K_8, \ \cdots, \ K_8  \, ) ,  \qquad
 K_{8} = {\rm diag}(  \, 1, \ -1, \ \cdots,  \ -1  \, ), \qquad
 t = m, \  n. \eqno (4.9)
$$  

It is easy to see from Eqs.(4.4) and (4.5) that the two kinds of block
 Kronecker products are actually  constructed by replacing all elements in
 the standard Kronecker product of matrices with $ 8 \times 8 $ matrices.
 Hence the operation properties on these  two  kinds of products are much
 similar to those on the standard Kronecker product of  matrices. We do not
  intend to list them here. 

\medskip

We next present some operation properties on the two real matrix
representations of octonion matrices.

\medskip

\noindent {\bf Theorem 4.1.} \,  {\em  Let $A,  \, B  \in \Bbb O^{ m \times
n}, \, \lambda \in \Bbb R$ be given.  Then

{\rm (a)} \ $ A=B  \Longleftrightarrow  \omega(A) = \omega(B) 
\Longleftrightarrow \ \nu(A) = \nu(B),$ i. e.$,$ $ \omega $ and $\nu$ 
are 1-1. 

{\rm (b)} \ $ \omega(A + B )=  \omega(A) + \omega(B), \ \ and \ \ 
\nu(A + B )= \nu(A) + \nu(B).$

{\rm (c)} \  $\omega( \lambda A) =  \lambda\omega(A), \ \ and \ \
\nu( \lambda A)= \lambda \nu(A).$

{\rm (d)}  \ $\omega (I_{m}) =I_{8m}, \ \ and \ \   \nu(I_{m}) = I_{8m}.$  

{\rm (e)}  \ $ \omega(A^*) = \omega^T(A), \ \ and \ \   \nu(A^*) = \nu^T(A),$
where $ A^* =  (\overline{a_{ts}}) $ is the conjugate transpose  of $ A$. }

\medskip

\noindent {\bf Theorem 4.2.} \,  {\em  Let $A \in \Bbb O^{ m \times n}$ be
given. Then
$$
A =  \frac{1}{8} E_{8m} \omega(A)E^T_{8n},    \eqno (4.10)
$$
where
$$ 
E_{8t} = {\rm diag}( \, E_8, \ \cdots, \ E_8  \,), \ \  and  \ \
 E_{8} = {\rm diag}( \, 1, \ e_1, \ \cdots,  \ e_7  \,), \ \ \  t = m, \ n.
$$  }
\noindent {\bf Proof.}  \ Follows  directly from Corollary 2.7.
 \qquad $ \Box$ 

\medskip

Since the multiplication of matrices over $\Bbb O$ is completely not
associative,  no identities on products of octonions matrices  can be
established over $\Bbb O$ in general. Consequently, no identities on
products of the two kinds of real matrix representations of octonion
matrices can be established.  In spit of this, we can still
 apply Eqs.(4.1) and (4.2) to deal with various problems related to octonion
matrices. Next are some results on the relationship of $\omega (\cdot)$, 
 $\nu (\cdot)$ and vec$(\cdot)$ for matrices over $ \Bbb O.$ 
     
\medskip

\noindent {\bf Lemma 4.3.} \, {\em  Let  $A \in \Bbb O^{n \times 1}, \, B
\in \Bbb O^{1 \times n}$ and  $ x \in \Bbb O$ be given. Then
$$
{\rm vec\,}(Ax) = \omega( A ) \overrightarrow{x}   \ \ and \ \
 {\rm vec\,}(xB) = \nu( B^T ) \overrightarrow{x}. \eqno (4.11)
$$ } 
\noindent {\bf Proof.}  \,  Let $ A = [\, a_1, \ \cdots , \ a_n \,]^T$ and
$ B = [ \, b_1, \ \cdots , \ b_n \, ]^T $. Then by Eqs.(2.3), (2.7) and
Eqs.(4.1)---(4.3) we find
 $$
{\rm vec}(Ax) = \left[ \begin{array}{c} \overrightarrow{a_1x}   \\ \vdots \\ \overrightarrow{a_nx} \end{array} \right]
= \left[ \begin{array}{c} \omega(a_1) \overrightarrow{x}   \\ \vdots \\  \omega(a_n) \overrightarrow{x} \end{array} \right]
= \left[ \begin{array}{c} \omega(a_1)  \\ \vdots \\  \omega(a_n) \end{array} \right]\overrightarrow{x} = 
\omega(A) \overrightarrow{x}, 
$$
and
$$
{\rm vec}(xB) = \left[ \begin{array}{c} \overrightarrow{xb_1}   \\ \vdots \\ \overrightarrow{xb_n} \end{array} \right]
= \left[ \begin{array}{c} \nu(b_1) \overrightarrow{x}   \\ \vdots \\  \nu(b_n) \overrightarrow{x} \end{array} \right]= \left[ \begin{array}{c} \nu(b_1)  \\ \vdots \\  \nu(b_n) \end{array} \right]\overrightarrow{x} = 
\nu(B^T) \overrightarrow{x}.  \qquad \Box
$$
\noindent {\bf Lemma 4.4.}  \, {\em  Let  $A \in \Bbb O^{m \times n}, \ X
\in \Bbb O^{n \times 1}$ and $a \in \Bbb O$ be given. Then
$$
{\rm vec}(AX) = \omega(A){\rm vec}X \ \ and \ \  {\rm vec}(Xa)
= [ \, \nu(a)\widehat{\otimes} I_{8n} \,]{\rm vec}X
= \nu(a)\widehat{\otimes}{\rm vec}X.  \eqno (4.12)
 $$ } 
\noindent {\bf Proof.} \ Let $ A = [\, A_1, \ \cdots , \ A_n \,]$ and
$ X = [\, x_1, \ \cdots , \ x_n \,]^T $. Then
 by Eq.(4.11) we find 
\begin{eqnarray*}
{\rm vec}(AX) & = & {\rm vec}( A_1x_1 + \cdots +  A_nx_n) \\
 & = & {\rm vec}( A_1x_1 ) + \cdots +  {\rm vec}( A_nx_n) \\
& = & \omega( A_1 ) {\rm vec}x_1 + \cdots +  \omega( A_n ) {\rm vec}x_n \\
& = & [\, \omega( A_1 ), \ \cdots, \ \omega( A_n ) \,]
\left[ \begin{array}{c} {\rm vec\,}x_1   \\ \vdots \\ {\rm vec\,}x_n
\end{array} \right] = \omega( A ) {\rm vec}X,
\end{eqnarray*}
as required for the first equality in (4.12). On the other hand,
$$
{\rm vec}(Xa) = \left[ \begin{array}{c} \overrightarrow{x_1a}   \\ \vdots
 \\  \overrightarrow{x_na} \end{array} \right]
 = \left[ \begin{array}{c} \nu(a) \overrightarrow{x_1}   \\ \vdots
 \\ \nu(a) \overrightarrow{x_n} \end{array} \right]
 = [ \, \nu(a)\widehat{\otimes} I_{8n} \,]{\rm vec}X
= \nu(a)\widehat{\otimes}{\rm vec}X,
$$
as required for the second equality in (4.12). \qquad $\Box$

\medskip

\noindent {\bf Lemma 4.5.} \, {\em  Let  $B \in \Bbb O^{p \times 1}$ and $X
\in \Bbb O^{n \times p}$ be given. Then
$$
{\rm vec\,}(XB) = [\, \nu(B^T) \widehat{\otimes}I_{8n} \,] {\rm vec}X . \eqno (4.13)
 $$ } 
\noindent {\bf Proof.}  \  Let $ X = [\, X_1, \ \cdots , \ X_p \,]$ and
$ B = [ \,  b_1, \ \cdots , \ b_p \,]^T $. Then it follows from the second
equality in (4.12) that
\begin{eqnarray*}
{\rm vec}(XB) &= & {\rm vec}( X_1b_1 + \cdots +  X_pb_p) \\
 & = & {\rm vec}( X_1b_1 ) + \cdots +  {\rm vec}( X_pb_p) \\
& = & ( \nu( b_1 )\widehat {\otimes} I_{8n} ) {\rm vec}X_1 + \cdots +
 ( \nu( b_p )\widehat {\otimes} I_{8n}){\rm vec}X_p \\
& = & ( \, [ \,  \nu( b_1 ), \ \cdots, \  \nu( b_p ) \, ]\widehat{\otimes}I_{8n}
\, ) \left[ \begin{array}{c} {\rm vec}X_1   \\ \vdots \\ {\rm vec}X_p
\end{array} \right] =  [ \nu( B^T )\widehat{\otimes}I_{8n} ]{\rm vec}X, 
\end{eqnarray*} 
as required for Eq.(4.13). \qquad  $\Box$ 

\medskip

Based on the above several lemmas, we can find the following three general
results. 

\medskip

\noindent {\bf Theorem 4.6.} \, {\em  Let  $A = ( a_{st}) \in \Bbb O^{m
\times n}$ and $  X \in \Bbb O^{n \times p}$ be given. Then
$$
{\rm vec}(AX) = [ \, I_{8p} \widehat{ \otimes} \omega( A ) \, ]{\rm vec}X.
\eqno (4.14)
$$  }
\noindent {\bf Proof.}  \,  Let $ X = [ \, X_1, \ \cdots , \ X_p \, ]$. Then
we find by Eq.(4.12) that
\begin{eqnarray*}
{\rm vec}(AX) &= & {\rm vec}[\, AX_1, \  \cdots, \   AX_p \,] \\
& = &  \left[\, {\rm vec}( AX_1), \  \cdots, \  {\rm vec}(AX_p)  \, \right]
\\
& = &  \left[ \,  \omega( A ){\rm vec}X_1, \  \cdots, \ \omega( A )
{\rm vec}X_p  \, \right] \\
& = & {\rm diag}( \, \omega( A ), \  \cdots, \  \omega( A ) \, )
[\, {\rm vec}X_1, \  \cdots, \ {\rm vec}X_p  \, ]
= [ \, I_{8p} \widehat{\otimes} \omega( A ) \, ] {\rm vec}X,
\end{eqnarray*}
establishing Eq.(4.14). \qquad  $\Box$

\medskip

\noindent {\bf Theorem 4.7.}  \, {\em  Let  $B = ( b_{st}) \in \Bbb O^{p
\times q}$ and $  X \in \Bbb O^{n \times p}$ be given. Then
$$
{\rm vec}(XB) = [\, \nu(B^T) \widehat{ \otimes} I_{8n} \,]{\rm vec}X.
\eqno (4.15)
$$ }
\noindent {\bf Proof.} \,  Let $ B = [ \, B_1, \ \cdots , \ B_q \, ]$.
Then we find by Eq.(4.13) that
\begin{eqnarray*}
{\rm vec}(XB) & = & {\rm vec}[\, XB_1, \  \cdots, \   XB_q \,] \\
 & = & \left[ \begin{array}{c} {\rm vec}XB_1   \\ \vdots  \\ 
{\rm vec}XB_q \end{array} \right] \\
 &= & \left[ \begin{array}{c}
 \left[ \, \nu (B^T_1) \widehat{\otimes}I_{8n} \, \right]{\rm vec}X  \\ \vdots \\
 \left[ \, \nu (B^T_q) \widehat{\otimes}I_{8n} \, \right] {\rm vec}X
  \end{array} \right] = \left[ \begin{array}{c} \left[ \, \nu (B^T_1) \widehat{\otimes}I_{8n} \, \right]  \\ 
\vdots
\\ \left[ \, \nu (B^T_q) \widehat{\otimes}I_{8n} \, \right]
\end{array} \right]{\rm vec}X =[ \, \nu (B^T) \widehat{ \otimes} I_{8n} \, ]{\rm vec}X,
\end{eqnarray*} 
as rerquired for Eq.(4.15). \qquad  $ \Box$

\medskip

\noindent {\bf Theorem 4.8.}  \, {\em  Let $A = ( a_{st}) \in \Bbb O^{m
\times n}, \ B = ( b_{st}) \in \Bbb O^{p \times q},$ and $  X \in
\Bbb O^{n \times p}$ be given. Then
$$
{\rm vec}[(AX)B] = [\, \nu (B^T) \widehat{ \otimes} \omega ( A ) \, ]
{\rm vec}X, \ \  and  \ \
 {\rm vec}[A(XB)] = [ \,  \omega( A ) \tilde { \otimes} \nu ( B^T ) \,
]{\rm vec}X. \eqno (4.16)
$$ }
\noindent {\bf Proof.} \,  According to Eqs.(4.14) and (4.15), we find
that
\begin{eqnarray*} 
{\rm vec}[(AX)B] & = & [ \, \nu(B^T) \widehat{ \otimes} I_{8m} \, ]
{\rm vec}(AX) \\
& = & [ \, \nu(B^T) \widehat{ \otimes} I_{8m} \, ][ \, I_{8p}
\widehat{ \otimes} \omega( A ) \, ]{\rm vec}X =
[ \, \nu(B^T) \widehat{ \otimes} \omega( A ) \, ]{\rm vec}X,
\end{eqnarray*} 
and
\begin{eqnarray*}
{\rm vec}[A(XB)] & = & [ \, I_{8p} \widehat{ \otimes} \omega( A ) \, ]
{\rm vec}(XB) \\
& = & [ \, I_{8p} \widehat{ \otimes} \omega( A ) \, ][ \, \nu(B^T)
\widehat{ \otimes} I_{8n}  \, ]{\rm vec}X = [ \,  \omega( A )
\tilde{ \otimes} \nu( B^T ) \, ]{\rm vec}X,
\end{eqnarray*} 
as required for Eq.(4.16). \qquad $\Box$

\medskip

\noindent {\bf Theorem 4.9.} \, {\em  Let $A = ( a_{st}) \in \Bbb O^{n
\times n}, \ X = ( b_{st}) \in \Bbb O^{n \times p}, \ Y \in \Bbb O^{q
\times n}$ be given$,$ and denote
$$ 
A^{(k|} * X  = A(A \cdots (AX) \cdots )), \ \ and \ \ 
  Y * A^{|k)} = (( \cdots (YA) \cdots ) A)A.
$$
Then
$$
{\rm vec}(A^{(k|} * X) = [\, I_{8p} \widehat{ \otimes} \omega^k( A ) \,
]{\rm vec}X,  \ \ and \ \ {\rm vec}( Y * A^{|k)} ) = [ \, \nu^k(A^T)
\widehat{ \otimes} I_{8q}  \,]{\rm vec}Y.  \eqno (4.17)
$$
}
Just as the standard Kronecker products for matrices over any field, the
three formulas in Eqs.(4.14)---(4.16) can  directly be used for transforming
any linear matrix equations over $ \Bbb O $ into an ordinary linear system
 of equation over  $ \Bbb R$. For example,
$$ \displaylines{
\hspace*{2cm}
 AX = B \ \Longleftrightarrow \  [  \, I \widehat{ \otimes} \omega( A )  \,
 ]{\rm vec} X = {\rm vec}B,  \hfill 
\cr
\hspace*{2cm}
XA = B \ \Longleftrightarrow \  [  \, \nu(A^T)  \widehat{ \otimes} I \, ]
{\rm vec} X = {\rm vec}B,     \hfill
 \cr
\hspace*{2cm}
 A(BX) = C \ \Longleftrightarrow \  [  \, I \widehat{ \otimes} \omega( A )
 \omega( B ) \, ]{\rm vec} X = {\rm vec}C, \hfill
 \cr
\hspace*{2cm}
(XA)B = C \ \Longleftrightarrow \  [  \, \nu(B^T)\nu(A^T) \widehat{ \otimes}
 I \, ] {\rm vec} X = {\rm vec}C, \hfill
 \cr
\hspace*{2cm}
(AX)B = C   \ \Longleftrightarrow \ [ \, \nu(B^T) \widehat{ \otimes}
\omega( A ) \, ]{\rm vec}X = {\rm vec}C, \hfill
 \cr
\hspace*{2cm}
 A(XB) = C  \ \Longleftrightarrow \  [\, \omega( A ) \tilde{ \otimes}
 \nu( B ) \, ]{\rm vec}X = {\rm vec}C,  \hfill
 \cr
\hspace*{2cm}
AX- XB = C \ \Longleftrightarrow \  [\, I \widehat{ \otimes} \omega( A )  -
\nu(B)\widehat{ \otimes} I \,] {\rm vec}X  = {\rm vec}C , \hfill
 \cr
\hspace*{2cm}
 (AX)A - A(XA) = B \ \Longleftrightarrow \ [ \, \nu(A^T)\widehat{ \otimes}
 \omega( A ) -  \omega( A ) \tilde{ \otimes} \nu( A^T ) \, ]  {\rm vec}X =
 {\rm vec}B.  \hfill
 \cr }
$$
Theoreticlly speaking, various problems related to linear matrix equations
over the octonion algebra now have a complete resolution. 

\medskip

Below are several simple results related to solutions of linear
 matrix equations over $\Bbb O$. 

 \medskip

\noindent {\bf Definition 4.3.} \,  Let $ A \in \Bbb O^{n\times n}$ be
given. If its left adjoint matrix $\omega( A )$ is invertible, then $ A $
is said to be {\em completely invertible}. 

\medskip

\noindent {\bf Theorem 4.10.} \, {\em  Let $A = ( a_{st}) \in \Bbb O^{m
\times m}$ and  $  B = ( b_{st}) \in \Bbb O^{m \times n}$ be given. If
$ A $ is completely invertible$,$ then the matrix equation
$$
 AX = B,  \eqno (4.18)
$$
has a unique solution over  $ \Bbb O$. In that case$,$ if the real
characteristic polynomial of  $\omega( A )$ is
$$
p(\lambda) =\lambda^t + r_{t-1}\lambda^{t-1} + \cdots +  r_{1}\lambda + r_0, \eqno (4.19)  
$$
where $r_0$ is the determinant of $ \omega(A),$ then the  unique solution of Eq.{\rm (4.18)} can be expressed as 
$$ 
X = - \frac{1}{r_0}[ \, A^{(t-1|}*B +  r_{t-1}( A^{(t-2|}*B ) +  \cdots +
 r_{3}A(AB) + r_{2}AB +  r_{1}B \,]. \eqno (4.20)
$$ }
\noindent {\bf Proof.} \ According to Eq.(4.14), the matrix equation (4.18)
is equivalent to
$$
 [\, I_{8n} \widehat{ \otimes} \omega( A )\,]{\rm vec} X = {\rm vec}B.
 \eqno (4.21)
$$
Because $ \omega( A )$ is invertible,  $ I_{8m} \hat{ \otimes} \omega( A)$
is also invertible. Hence the solution of Eq.(4.25) is unique and this
solution is
$$  
{\rm vec} X = [\, I_{8n} \widehat{ \otimes} \omega( A ) \,]^{-1}{\rm vec}B
= [\,I_{8m} \widehat{ \otimes} \omega^{-1}( A ) \,]{\rm vec}B.
$$ 
Observe that 
$$
\omega^t( A ) + r_{t-1}\omega^{t-1}( A ) + \cdots +  r_1 \omega( A )  +
r_0I_{8m} = 0
$$ 
holds. We then have 
$$
\omega^{-1}( A ) = - \frac{1}{r_0} \left[ \, \omega^{t-1}( A ) +
r_{t-1}\omega^{t-2}( A ) + \cdots +  r_2 \omega( A )  +
r_1 I_{8m} \, \right].
$$ 
Thus
$$
I_{8n} \widehat{\otimes} \omega^{-1}( A ) =
 - \frac{1}{r_0}[ \, I_{8n} \widehat{ \otimes}\omega^{t-1}( A ) +
 r_{t-1}(\, I_{8n} \widehat{ \otimes}\omega^{t-2}( A ) \, ) + \cdots +
 r_2 ( \, I_{8n} \widehat{ \otimes}\omega( A ) \,)  +  (\, r_1 I_{8n}
 \widehat{ \otimes }I_{8m} \,) \,],
$$
and 
\begin{eqnarray*}  
{\rm vec \,} X  =   [\, I_{8n} \widehat{\otimes} \omega( A ) \,]^{-1}{\rm
vec\,}B & = & - \frac{1}{r_0}[ \, ( \, I_{8n} \widehat{ \otimes}\omega^{t-1}( A )
 \, ) {\rm vec\,}B  + r_{t-1}(\, I_{8n} \widehat{ \otimes}\omega^{t-2}( A ) \,
  ){\rm vec}B \\
 &  & \ + \cdots +  r_2 ( \, I_{8n} \widehat{ \otimes}\omega( A ) \,)
 {\rm vec\,}B + r_1 (\, I_{8n} \widehat{ \otimes }I_{8m} \,){\rm vec\,}B \,].
\end{eqnarray*} 
Retuning it to octonion matrix expression by Eq.(4.17), we obtain Eq.(4.24).
   \qquad  $ \Box $  

\medskip

Similarly we have the following. 

\medskip

\noindent {\bf Theorem 4.11.} \, {\em  Let $A = ( a_{st}) \in \Bbb O^{m \times m}$ and $ B = ( b_{st}) \in \Bbb O^{n \times m}$ be given. If  $ A $ is completely invertible, then the matrix equation $ XA= B$ has a unique solution over  $ \Bbb O$. In that case, if the real characteristic polynomial of $\omega( A )$ is
$$ 
p(\lambda) =\lambda^t + r_{t-1}\lambda^{t-1} + \cdots +  r_{1}\lambda + r_0, \eqno (4.22) 
$$
then the  unique solution of $ XB = A $ can be expressed as 
$$ 
X = - \frac{1}{r_0}[ \, B* A^{|t-1)} +  r_{t-1}( B * A^{|t-2)} ) +  \cdots +
  r_3(BA)A + r_{2}BA +  r_{1}B \,]. \eqno (4.23)
$$ } 

For simplicity, the two solutions in Eqs.(4.20) and (4.23) can also be
written as
$$ 
X = L^{-1}_A\circ B,  \qquad X = B\circ R^{-1}_A, \eqno (4.24)
$$ 
where $ L^{-1}_A $ and  $R^{-1}_A$ are, respectively, called the
{\em left and the right inverse operators} of the completely invertible
octonion matrix $ A $. Some properties on these two inverse operators
are listed below.

\medskip

\noindent {\bf Theorem 4.12.} \, {\em  Let $A  \in \Bbb O^{m \times m}$ be
an  completely invertible matrix$,$ $B \in \Bbb O^{m \times n}$ and
$ C \in \Bbb O^{n \times m}$ be given. Then
$$ 
A( L^{-1}_A \circ B) = B,   \qquad  A( L^{-1}_A \circ I_m) = I_m,
\eqno (4.25)
$$ 
$$ 
 L^{-1}_A \circ (AB) = B,  \qquad  L^{-1}_A \circ A = I_m, \eqno (4.26)
$$ 
$$ 
 ( C \circ R^{-1}_A )A = C,  \qquad  (  I_m \circ R^{-1}_A ) A  = I_m, \eqno (4.27)
$$ 
$$ 
 (CA)\circ R^{-1}_A = C,  \qquad  A \circ R^{-1}_A  = I_m. \eqno (4.28)
$$ }
\noindent {\bf Proof.} \, Follows from Theorems 4.10 and 4.11. \qquad
 $ \Box $

\medskip

We can also consider the inverses of  octonion matrices in the usual sense.
Let  $A \in \Bbb O^{m \times m}$ be given. If there are $ X, \, Y \in
\Bbb O^{m \times m}$ such that $ XA = I_m $ and $ AY = I_m $, then $ X $
and $ Y $ are, respectively, called the {\em left inverse} and the 
{\em right inverse} of $ A $, and denoted by $  A_L^{-1} := X$ and  $ A_R^{-1} := Y$. From Theorems
4.10 and 4.11, we know that a square matrix of order $ m $ over
$ \Bbb O$ has a left inverse  if and only if the equation $ [ \, \nu(A^T)
\widehat{ \otimes}I_{8m} \,] {\rm vec} X = {\rm vec}I_m $ is
solvable, and $ A $ has a right inverse  if and only if
the equation $ [ \, I_{8m} \widehat{ \otimes} \omega(A) \,]
{\rm vec}Y = {\rm vec}I_m $ is solvable. These two facts imply that the left
and the right inverses of a square matrix may not be  unique, even  both of
them exist. As two special cases, we have the following.

\medskip

\noindent {\bf Theorem 4.13.} \, {\em  Let $A \in \Bbb O^{m \times m}$ be
given.  Then the left and the right inverses  of $ A $ are unique if and only
 if  $ A $   is completely invertible. In that case$,$ if the real
 characteristic polynomial of  $\omega( A )$ is
$$
p(\lambda) =\lambda^t + r_{t-1}\lambda^{t-1} + \cdots +  r_{1}\lambda + r_0, 
$$
then the unique left and the unique right inverses $ A $  can be expressed
as
$$ 
A_L^{-1} = - \frac{1}{r_0}[ \, A^{(t-1|} +  r_{t-1} A^{(t-2|}  +  \cdots +
r_{3}A(A^2) + r_{2}A^2 + r_{1}I_m \,],
$$  
and
$$
A_R^{-1} = - \frac{1}{r_0}[ \, A^{|t-1)} +  r_{t-1}A^{|t-2)} +  \cdots + r_{3}(A^2)A + r_{2}A^2 +  
r_{1}I_m \,],
$$  
where $ A^{(s|} := A(A( \cdots (AA) \cdots ))$ and $  A^{|s)} :=
(( \cdots (AA) \cdots )A)A.$  } 

\medskip

\noindent {\bf Proof.} \ Follows directly from Theorems 4.10 and 4.11.
 \qquad $ \Box $

\medskip

Based on Theorems 4.10 and 4.12, as well as Eqs.(4.25)---(4.28), we can
also derive the following two simple
results. 

\medskip

\noindent {\bf Corollary  4.14.} \, {\em If $A \in \Bbb O^{m
\times m}$ is completely invertible$,$ and $ AB_1 = AC_1$ and $ B_2A = C_2A,$
then $ B_1 = C_1 $ and $ B_2 = C_2.$ In other words$,$ the left and the right
cancellation rules hold for completely invertible  matrices. }

\medskip

\noindent {\bf Corollary 4.15.}  \,  {\em  Suppose that $A \in \Bbb O^{m
\times m}, \  B \in \Bbb O^{n \times n}$
 are completely invertible and $ C \in \Bbb O^{m \times n}$. Then

 {\rm (a)} \ The  matrix equation $ A(XB) = C $ has a unique solution
 $ X = ( L^{-1}_A \circ C )R^{-1}_B.$

 {\rm (b)} \  The matrix equation $ (AX)B = C $ has a unique solution
  $ X = L^{-1}_A ( C \circ R^{-1}_B ), $ \\
where $ L^{-1}_A $ and  $ R^{-1}_B$ are the left and the right inverse
operators of $ A $ and $ B $ respectively. }

\medskip

Our next result is concerned with the extension of the Cayley-Hamilton
theorem to octonion matrices, which could be regarded as one of most
successful applications of  matrix representations of octonions.

\medskip

\noindent {\bf Theorem 4.16.} \, {\em  Let $A \in \Bbb O^{m \times m}$ be
given and suppose that the real characteristic polynomial of $\omega( A )$
is
$$
p(\lambda) =\lambda^t + r_{t-1}\lambda^{t-1} + \cdots +  r_{1}\lambda + r_0. 
$$
Then $ A $ satisfies the following two identities
$$
 A^{(t|} +  r_{t-1}A^{(t-1|} +  \cdots + r_{3} A(AA)  + r_{2}A^2 +  r_1A  +
 r_0I_m = 0,  \eqno (4.29)
$$ 
$$
 A^{|t)} +  r_{t-1}A^{|t-1)} +  \cdots+ r_{3}(AA)A  + r_{2}A^2 +  r_1A  +
 r_0I_m = 0.  \eqno (4.30)
$$ }  
\noindent {\bf Proof.} \, Observe that $ p[\omega( A )] = 0$. It follows
that
 $$ 
[ \, I_{8m} \widehat{ \otimes}  p[\omega( A ) ] \,] {\rm vec}I_m = 0.  \eqno (4.31)
$$
On the other hand, it is east to see by Eq.(4.17) that
$$
 {\rm vec} A^{(s|} = {\rm vec}(  A^{(s|} * I_m ) =
 [\, I_{8m} \widehat{\otimes} \omega^s( A ) \,]
{\rm vec}I_m, \ \ \ s = 1, \, 2, \,  \cdots.
$$
Thus we find that
$$
\displaylines{ 
\hspace*{0cm}
[ \, I_{8m} \widehat{ \otimes} p(\omega( A ) ) \,] {\rm vec\,}I_m  \hfill
\cr 
\hspace*{0cm}  
= [\, I_{8m} \widehat{ \otimes} \omega^t( A ) + r_{t-1}(\, I_{8m}
\widehat{ \otimes} \omega^{t-1}( A ) \,) + \cdots +  r_{1}(\, I_{8m}
\widehat{ \otimes} \omega( A ) \,) +  \ r_0(\, I_{8m} \widehat{ \otimes}I_{8m} \,) \,] {\rm vec\,}I_m  \hfill
\cr 
\hspace*{0cm} 
= (\, I_{8m} \widehat{ \otimes} \omega^t( A ) \,){\rm vec}I_m +
r_{t-1}(\, I_{8m} \widehat{ \otimes} \omega^{t-1}( A ) \,) {\rm vec\,}I_m  +
\cdots +  r_{1}(\, I_{8m} \widehat{ \otimes} \omega( A ) \,){\rm vec\,}I_m
 + r_0(\, I_{8m} \widehat{ \otimes}I_{8m} \,) {\rm vec\,}I_m  \hfill
\cr 
\hspace*{0cm} 
= {\rm vec \,} A^{(t|} +  r_{t-1}{\rm vec\,}A^{(t-1|} +  \cdots +
 r_1{\rm vec\,}A  + r_0{\rm vec\,}I_m  \hfill
\cr 
\hspace*{0cm} 
= {\rm vec} [\,  A^{(t|} +  r_{t-1}A^{(t-1|} +  \cdots +  r_1A  +
 r_0I_m \,]. \hfill
}
$$
The combination of this equality with Eq.(4.31) results in Eq.(4.29).
The identity in Eq.(3.30) can be established similarly.   \qquad  $ \Box $

\medskip

Finally we present a result on real eigenvalues of Hermitian
octonion matrices. 

\medskip

\noindent {\bf Theorem 4.17.} \, {\em Suppose that $ A \in \Bbb O^{ m \times m }$ is Hermitian$,$ that is$,$  $ A^* = A$. Then $A$ and its real adjoint 
 $\omega(A)$ have identical real eigenvalues. } 

\medskip

\noindent { \bf  Proof.} \, Since $ A = A^* $, we know by Theorem 4.1(e) that 
 $\omega (A) = \omega (A^*) = \omega^T(A),$ that is, $ \omega (A)$ is a real 
symmetric matrix. In  that case, all eigenvalues of $ \omega (A)$ are real.  Now suppose that 
$$
\omega (A)X = X \lambda, \eqno (4.32)
$$   
where $ \lambda \in \Bbb R$ and $ X \in \Bbb R^{ 8m \times 1}$. Then there is
 unique  $ Y \in \Bbb O^{ m \times 1}$ such that 
${\rm vec \,}Y = X$. In that case, it is easy to find  by Theorem 4.1(a) and 
Eq.(4.12) that 
$$
\omega (A)X = X \lambda  \Longrightarrow  
\omega (A){\rm vec \,}Y = {\rm  vec \,}Y \lambda \Longrightarrow 
{\rm vec \,}(AY) = {\rm vec \,}(Y \lambda ) \Longrightarrow  AY = Y \lambda, \eqno (4.33)
$$
which implies that $ \lambda $ is a real eigenvalue of $ A $, and $ Y $ is a eigenvector of $ A $ corresponding  to this $ \lambda$. Conversely  suppose that 
 $ AY = Y \lambda,$ where $ \lambda \in \Bbb R$,  $ Y \in \Bbb O^{m \times 1 }$. Then taking vec operation on its both sides according to Eq.(4.12) yields
$$
\omega (A) {\rm vec \,} Y  = {\rm vec \,}Y \lambda.
$$
This implies that $\lambda $ is also a real eigenvalue of $ \omega (A)$ and 
${ \rm vec \, } Y $ is a real eigenvector of $ \omega ( A)$ associcated with this $\lambda$.  \qquad $ \Box$

\medskip

The above result clearly shows that real eigenvalues and  the corresponding eigenvectors of 
a Hermitian octonion matrix $ A $ can all be determined by its real adjoint $\omega (A)$. 
 Since $\omega (A) $ is a real symmetric $ 8m \times 8m$  matrix, it has $ 8m $ eigenvalues 
and $8m$ corresponding orthogonal eigenvectors.  

Now a fundamental problem would naturally be asked: how many different real eigenvalues can  a Hermitian octonion matrix $ A $ have at most? For a 
$ 2 \times 2 $ Hermitian  octonion matrix 
$A = \left[\begin{array}{cc} a & b  \\ \overline{b} & c  \end{array} \right]$, 
where $ a, \ c \in  \Bbb R$,  its real adjoint is 
$$
\omega (A) = \left[\begin{array}{cc} aI_8 & \omega (b)  \\ \omega^T(b) & cI_8  
\end{array} \right]. 
$$
Clearly the characteristic polynomial of $ \omega (A)$ is 
$$ 
{\rm det}( \, \lambda I_{16} - \omega (A) \, ) =  [ \, (\lambda - a)(\lambda - c) - |b|^2 \, ]^8. 
$$
This shows that $ \omega (A)$, and correspondingly $ A$, has 2 eigenvalues, each of which 
has a multiplicity 8.       

The eigenvalue problem for $ 3 \times 3 $ Hermitian octonion matrices was recently examined 
by Dray and Manogue \cite{DM} and Okubo \cite{Ok1}. They showed by algebraic methods that  
every $ 3 \times 3 $ Hermitian octonion matrix has 24 real eigenvalues which are divided into 6 groups, 
each of them has multiplicity 4. Now according to Theorem 4.17, the  real eigenvalues of any $ 3 \times 3 $ Hermitian 
octonion matrix  
$$
A = \left[\begin{array}{ccc} a_{11} & a_{12} & a_{13}  
\\ \overline{a}_{12} & a_{22} & a_{23} \\  \overline {a}_{13} & 
\overline {a}_{23} & a_{33}  \end{array} \right],   \qquad  a_{11}, \ a_{22}, \ 
a_{33} \in \Bbb R,
$$   
can be completely determined  by its real adjoint  
$$
\omega (A) = \left[\begin{array}{ccc} \omega (a_{11}) & \omega (a_{12}) & \omega (a_{13})  
\\ \omega^T({a}_{12}) & \omega (a_{22}) & \omega (a_{23}) \\  \omega^T({a}_{13}) & 
\omega^T({a}_{23}) & \omega (a_{33})  \end{array} \right]. 
$$
Obviously this matrix  has 24 real eigenvalues and the $24$ corresponding real orthogonal eigenvectors. 
Numerical computation shows that these 24 eigenvalues are divided into 6 groups, each of them has multiplicity 4, 
which is consistent with the fact revealed in \cite{DM} and \cite{Ok1}. Moreover the 24 
real orthogonal eigenvectors can also be converted to octonion expressions by (4.33).

Furthermore, numerical computation reveals an interesting fact that the 32 real eigenvalues any 
$ 4 \times 4 $ Hermitian octonion matrix are divided into 16 groups, each of them has multiplicity 
2; the 40 real eigenvalues of any $ 5 \times 5$ Hermitian octonion matrix are divided 
into 20 groups, each of them has multiplicity 2.   

In general, we guess that for any $ m \times m $ Hermitian octonion matrix with $ m > 3$, 
its $8m$ real eigenvalues can be divided into $4m$ groups, each of them has  multiplicity 2.   

As a subsequent work of Thereom 4.17, one might naturally ask hwo to establish a possible factorization for a Hermitian octonion matrix using its real 
eigenvalues and corresponding octonion orthogonal eigenvectors, speak more precisely, 
 for an $ m \times m $ Hermitian octonion matrix  $A$, how construct a complete invertible 
octonion matrix $P$ (unitary?) and a real diagonal matrix $D$ such that $ A = PDP^{-1}$ using 
its $8m$ real eigenvalues and $8m$ corresponding octonion orthogonal eigenvectors. However, this problem 
seems quite curious, because the number of different real eigenvalues of an
 Hermitian octonion matrix is more than its order. This problem is also quite challenging, because various  
 traditional methods in associative matrix theory are not applicable to this non-associative case.  
 
\medskip

As pointed out in \cite{DM}, Hermitian octonion matrices can also have non-real   
right eigenvalues. Theoretically speaking, the non-real eigenvalue problem of Hermitian octonion matrices may also be converted to a problem related to 
real representations of octonion matrices. In fact, suppose that $ AX = X \lambda,$ where $\lambda \in \Bbb O$ and $ X \in \Bbb O^{ m \times 1}$. Then according to Eq.(4.12), it is equivalent to 
$$
\omega (A) {\rm vec \,}X = \nu (\lambda) \widehat{ \otimes}{\rm vec \,}X,
$$
or alternatively
$$
[\, \omega (A) - {\rm diag}( \, \nu (\lambda), \ \cdots, \ \nu (\lambda) \, ]{\rm vec \,}X = 0.
$$  
How to find $ \nu (\lambda)$ satisfying the equation remains to further study.

\medskip
 
\noindent {\bf Conclusions.} \, In this paper, we have introduced two
pseudo real matrix representations for octonions. Based on them we have  made 
a complete investigation to their operation properties and have considered their
 various applications to octonions and matrices of octonions. However our
 work could only be regarded as a first step in the  research of octonion 
matrix analysis  and its applications. Numerous problems related  to matrices of 
octonions remain to further examine, such as:

\begin{itemize}

\item[(a)] How to determine eigenvalues and eigenvectors of a square octonion
matrix, not necessarily Hermitian,  and what is the relationship of eigenvalues and eigenvectors of
a octonion matrix and its real adjoint matrices? 

\item[(b)] Besides Eq.(4.29) and (4.30), how to establish some other identities for octonion matrices through their adjoint matrices? 

\item[(c)] How to establish similarity theory for octonion matrices,
 and  how to determine the relationship between the similarity of octonions
  matrices and the similarity of their adjoint matrices?

\item[(d)]   How to consider various possible decompositions of octonion matrices,
 such as, LU decomposition, singular value decomposition and Schur decomposition?

\item[(e)]  How to characterize various particular octonion matrices, such as,
     idempotent matrices, nipoltent matrices, involutary matrices, unitary matrices,  
normal matrices, and so on?

\item[(f)]   How to define generalized inverses of octonion matrices when they
are not completely invertible?
\end{itemize}
\noindent and so on. As mentioned in the beginning of the section,
matrix multiplication for octonion matrices is completely not associative.
In that case, any further research to problems related matrices of octonions
is extremely difficult, but is also quite challenging. Any advance in 
solving the problems mentioned above could lead to remarkable new development in 
the real octonion algebra and its applications in mathematical physics. 

\medskip

Finally  we should point out that the results obtained in the
paper can use to establish pseudo matrix representations for
 real sedenions, as well as, in general,  for elements in  any  
$2^n$-dimensional real Cayley-Dickson algebras. \\

\end{document}